\documentclass[11pt]{article}
\usepackage[margin=4.46cm]{geometry}
\usepackage{graphicx}
\usepackage{epsfig}
\usepackage{amsmath}
\usepackage{amssymb}
\usepackage{caption}
\newcommand{\commentout}[1]{}

\def \Rset {{\mathbb R}}
\def \Cset {{\mathbb C}}
\def \Zset {{\mathbb Z}}

\def \Nset {{\mathbb N}}

\newcommand{\nit}{\noindent}
\newcommand{\no}{\nonumber}
\newcommand{\be}{\begin{equation}}
\newcommand{\ee}{\end{equation}}
\newcommand{\ba}{\begin{eqnarray}}
\newcommand{\ea}{\end{eqnarray}}
\newcommand{\bi}{\begin{itemize}}
\newcommand{\ei}{\end{itemize}}
\newcommand{\br}{\begin{eqnarray}}
\newcommand{\er}{\end{eqnarray}}

\newcommand{\qed}{\mbox{$\square$}\newline}

\newtheorem{theo}{Theorem}[section]

\newtheorem{lem}{Lemma}[section]

\newtheorem{rmk}{Remark}[section]

\begin{document}
\title{Curvature effect in shear flow: 
slowdown of turbulent flame speeds with Markstein number}
\author{Jiancheng Lyu, \ Jack Xin, \ Yifeng Yu \thanks{Department of Mathematics, University of California at Irvine, Irvine, CA 92697.
Email: (jianchel,jack.xin,yifengy)@uci.edu. The work was partly supported by 
NSF grants DMS-1211179 (JX), DMS-0901460 (YY), and CAREER Award DMS-1151919 (YY).}}
\date{}
\maketitle

\begin{abstract}
It is well-known in the combustion community that curvature effect in general 
slows down flame propagation speeds because it smooths out 
wrinkled flames. However, such a folklore has never been justified rigorously. 
In this paper, as the first theoretical result in this direction, 
we prove that the turbulent flame speed (an effective burning velocity) is 
decreasing with respect to the curvature diffusivity (Markstein number) for shear flows 
in the well-known G-equation model.  Our proof involves several 
novel and rather sophisticated inequalities arising from the nonlinear structure of the equation. On a related fundamental issue, we solve the selection problem of weak solutions or find the ``physical fluctuations" 
when the Markstein number goes to zero and solutions 
approach those of the inviscid G-equation model. The limiting solution 
is given by a closed form analytical formula. 
\end{abstract}
\vspace{.2 in}

\hspace{.12 in} {\bf AMS Subject Classification:} 70H20, 76M50, 76M45, 76N20.
\bigskip

\hspace{.12 in} {\bf Key Words:} Flame speeds, curvature smoothing, shear flows, 

\hspace{.12 in} speed slow-down, 
zero curvature limit.

\thispagestyle{empty}
\newpage

\section{Introduction}
\setcounter{equation}{0}
\setcounter{page}{1}

The curvature effect in turbulent combustion was first studied by Markstein \cite{Mar}, which says that if the flame front bends toward the cold region (unburned area, point C in Figure 1 below), the flame propagation slows down. If the flame front bends toward the hot spot (burned area, point B in Figure 1), it burns faster.
\begin{center}
\includegraphics[scale=0.6]{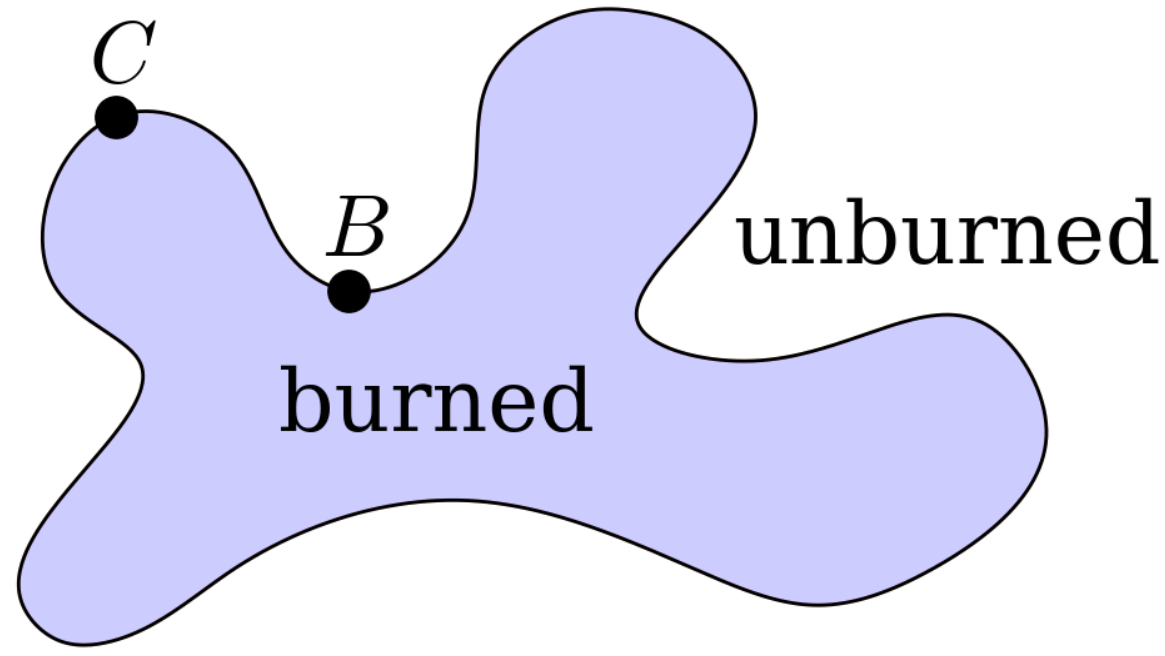}
\captionof{figure}{Curvature effect}
\end{center}

Below is an empirical linear relation  proposed by Markstein \cite{Mar} to approximate the dependence of the laminar flame speed  $s_l$  on the curvature (see also \cite{P2000}, \cite{Se1985}, etc):
\be\label{laminar}
s_l=s_{l}^{0}(1-\tilde d\; \kappa). 
\ee
Here $s_{l}^{0}$, the mean value, is a positive constant. The parameter 
$\tilde d>0$ is the so called Markstein length which is proportional to the flame thickness. 
The mean curvature along the flame front is $\kappa$. 

In general, $\kappa$ changes sign along a curved flame front. 
So a mathematically interesting and physically important question is:

\medskip 

\nit {\bf Q1:} {\it How does the ``averaged" flame propagation speed depend on the curvature term?}

\medskip

Of course, we first need to properly define an ``averaged speed", which is basically to average fluctuations caused by both the flow and the curvature. The theory of homogenization provides such a rigorous mathematical framework in environments with microscopic structures. In this paper, we employ the popular G-equation model in combustion community. 

Let the flame front be the zero level set of a reference function $G(x,t)$, where 
the burnt and unburnt regions are $\{G(x,t) < 0\}$ and $\{G(x,t) > 0\}$, respectively. 
See Figure 2 below. The velocity of ambient fluid $V:\Rset^n\to \Rset^n$ is 
assumed to be smooth, $\Zset^n$-periodic and incompressible (i.e. $div(V)=0$). 
The propagation of flame front obeys a simple motion law: $\vec {v}_{n}=s_l+V(x)\cdot n$, 
i.e., the normal velocity is the laminar flame speed ($s_l$) plus the projection of $V$ 
along the normal direction. This leads to the so--called $G$-equation, a level-set 
PDE \cite{OF2002,P2000}:
$$
G_t + V(x)\cdot DG + s_l \, |DG|=0 \quad \text{in $\Rset^n\times (0,+\infty)$}. 
$$

\begin{center}
\includegraphics[scale=0.8]{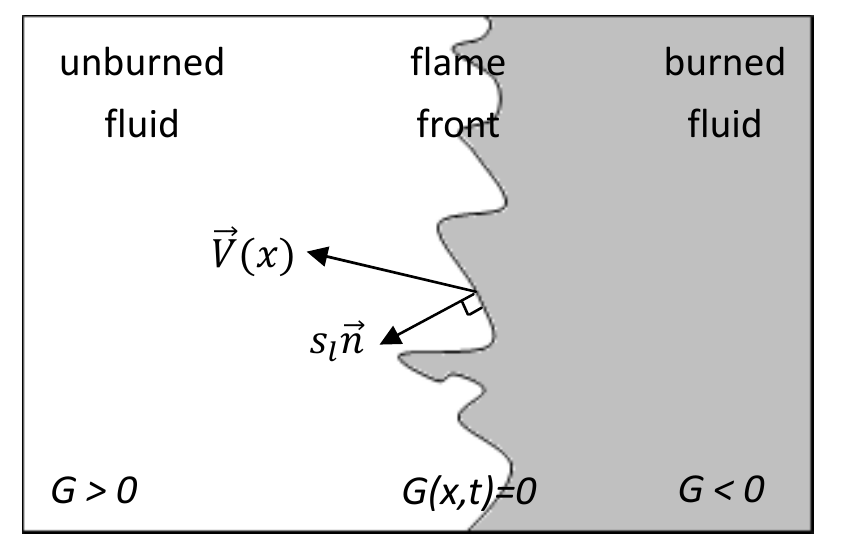} \quad \qquad
\captionof{figure}{Level-set formulation of front propagation}
\end{center}

\bigskip

Plugging the expression of the laminar flame speed (\ref{laminar}) 
into the G-equation and normalizing the constant $s_{l}^{0}=1$, we obtain a mean curvature type equation 
\be
G_t + V(x)\cdot DG + |DG|-\tilde d\, |DG|\, \mathrm{div}\left({DG\over |DG|}\right)=0.\label{ge1}
\ee

Turbulent combustion usually involves small scales. As a simplified model, 
we rescale $V$ as $V=V({x\over \epsilon})$ and write $\tilde d=d\epsilon$. 
Here $\epsilon$ denotes the Kolmogorov scale (the small scale in the flow). 
The diffusivity constant $d>0$ is called the Markstein number.  
We would like to point out that the dimensionless Markstein number is 
$d\cdot {\delta_L \over \epsilon}$ with $\delta_L$ denoting the flame thickness \cite{P2000}. 
In the thin reaction zone regime, $\delta_L = O(\epsilon)$, see Eq. (2.28) and Fig. 2.8 of \cite{P2000}. 
Without loss of generality, let ${\delta_L \over \epsilon}=1$. Then (\ref{ge1}) becomes
\be
G_{t}^{\epsilon} + V({x\over \epsilon})\cdot DG^{\epsilon} + |DG^{\epsilon}|-d\, \epsilon \, |DG^{\epsilon}|\;\mathrm{div}\left({DG^{\epsilon}\over |DG^{\epsilon}|}\right)=0.\label{ge1b}
\ee

Since $\epsilon \ll 1$, it is natural to look at $\lim_{\epsilon\to 0}G^{\epsilon}$, i.e., 
the homogenization limit. If for any $p\in \Rset^n$, 
there exists a unique number $\overline H_d(p)$ such that the following cell problem has (approximate) $\Zset^n$-periodic viscosity solutions in $\Rset^n$:
\be\label{curcell}
-d\, |p+Dw|\; \mathrm{div}\left({p+Dw\over |p+Dw|}\right)+ |p+Dw|+ V(y)\cdot (p+Dw)=\overline H_d(p),
\ee
then standard tools in the homogenization theory imply that 
$$
\lim_{\epsilon \to 0} G^{\epsilon}(x,t)=\bar G(x,t) \quad \text{locally uniformly in $\Rset\times [0, +\infty)$}.
$$
Here $\bar G$ is the unique solution to the following effective equation, 
which captures the propagation of the mean flame front (see Figure 3 below). 
\be
\begin{cases}
\bar G_{t}+\overline H_d(D\bar G)=0\\
\bar G(x,0)=G_0(x)  \quad \text{initial flame front}.
\end{cases}
\ee
\begin{center}
\includegraphics[scale=0.8]{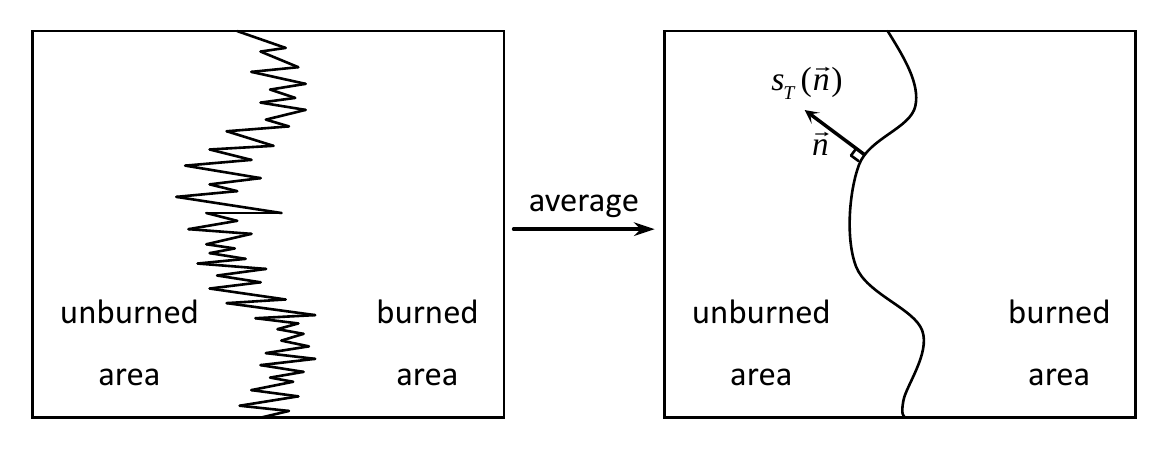}
\captionof{figure}{Average of fluctuations in the homogenization limit}
\end{center}

\bigskip

Solution to the cell problem (\ref{curcell}) formally describes fluctuations around 
the mean flame front, i.e., 
$$
G(x,t)=\bar G(x,t)+\epsilon w(x,{x\over \epsilon})+O(\epsilon^2),
$$
where for fixed location-time $(x,t)$ and $p=D\bar G(x,t)$, $w(x,\cdot)$ is a 
solution to (\ref{curcell}) with mean zero, i.e.,  $\int_{0}^{1}w(x,y)\,dy=0$.   
The quantity $\overline H_d(p)$, if it exists, can be viewed as  
the turbulent flame speed ($s_T(p)$) along a given direction $p$. 
There is a consensus in combustion literature that the 
curvature effect slows down flame propagation \cite{Ron}. 
Heuristically, this is because the curvature term smooths out the flame front and 
reduces the total area of chemical reaction \cite{Se1985}. However, 
this folklore has never been rigorously justified mathematically. 
If the curvature term is replaced by the full diffusion (i.e. the Laplacian $\Delta$), a 
dramatic slow-down is proved in \cite{LXY2011} for two dimensional cellular flows. 
So in the G-equation setting, {\bf Question 1} can be formulated as 

\medskip

\nit {\bf Q2}: {\it How does $\overline H_d(p)$ depend on the Markstein number $d$? In particular, is it decreasing with respect to $d$?}

\medskip

We remark that the decrease of turbulent flame speed with respect to 
the Markstein number has been experimentally observed (e.g., \cite{CWL2013}).

\subsection{Slow-down of Flame Propagation}

For general $V$, we do not even know the existence of $\overline H_d(p)$, i.e., the well-posedness of (\ref{curcell}). In fact, given the counter-example in \cite{CM2014} for a coercive mean curvature type equation, the cell problem (\ref{curcell}) and the homogenization in our non-coercive setting is very likely not well-posed in general. To avoid this existence issue, as the first step to investigate the above {\bf Question 2}, we consider the shear flow in this paper:
$$
V(x)=(v(x_2),0) \quad \text{for $x=(x_1,x_2)\in \Rset^2$.}
$$
Here $v:\Rset\to \Rset$ is a smooth periodic function. Then for $p=(\gamma, \mu)\in 
\Rset^2$, the cell problem (\ref{curcell}) is reduced to the following ODE:
\be\label{1d-curvature}
-{d\gamma^2w''\over \gamma^2+(\mu+w')^2}+\sqrt {\gamma^2+(\mu+w')^2}+\gamma v(y)=\overline H_d(p) \quad \text{in $\Rset$}.
\ee
It is then easy to show that there exists a unique number $\overline H_d(p)$ such 
that the ODE (\ref{1d-curvature}) has a  $C^2$ periodic 
solution. Throughout this paper,   we denote  $w$ as the unique solution satisfying that $w(0)=0$. To simplify notations, we omit the dependence of $w$ on $d$. The following is our main result.

\begin{theo}\label{main1} Assume that $v=v(y)$ is not a constant function. Then

(1) $\overline H_d(0, \pm \mu)=|\mu|$;

(2) {\bf (Major Part).} If $\gamma\not= 0$, 
$$
{\partial \overline H_d(p)\over \partial d}<0.
$$
So $\overline H_d$ is strictly decreasing with respect to the Markstein number $d$. 

(3) $\lim_{d\to 0^{+}}\overline H_d=\overline H_0$. Here $\overline H_0(p)$ is the unique number (effective Hamiltonian) such that the following inviscid equation admits periodic viscosity solutions
$$
\sqrt {\gamma^2+(\mu+w_{0}^{'})^2}+\gamma v(y)=\overline H_0(p)  \quad \text{in  $\Rset$}.
$$

(4) $\lim_{d\to +\infty}\overline H_d=|p|+\gamma \int_{0}^{1}v(y)\,dy$ and $\lim_{d\to +\infty}w=0$ uniformly in $\Rset$.
\end{theo}

Proofs for (1), (3) and (4) are simple. The real challenge 
is to prove the major part (2). A key step in our proof is to establish a 
highly sophisticated class of inequalities, see Lemma \ref{discrete} (the discrete version) 
and Theorem \ref{continuous} (a specific continuous version). 
Some calculations in high dimensions will be presented in Section 2.2 
when the ambient fluid is near rest. 

It might be tempting to think that there exists an explicit 
formula of $\overline H_d(p)$ since (\ref{1d-curvature}) is ``just" an ODE. However, 
this is not the case. For example, let us look at a simpler cell problem associated with 
the 1-d viscous Hamilton-Jacobi equation arising from large deviations and quantum mechanics:
$$
-d\, w''+|p+w'|^2+G(y)=\overline H(p,d) \quad \text{in $\Rset$}.
$$
Here the potential $G$ is a smooth periodic function and $\overline H(p,d)$ is the 
unique number such that the above equation has $C^2$ solutions. The 
viscous effective Hamiltonian $\overline H(p,d)$ actually determines the 
spectrum of the 1-d Schr\"odinger operator ($Lu=-du''+Gu$) and it is 
closely related to the inverse scattering solution of the KdV equation \cite{LTY2016}.  We want to remark that the strict decreasing of $\overline H(p,d)$ with respect to $d$ can be easily established in any dimension. See (\ref{viscous-basic}) in Remark \ref{vanishing-viscosity}.

\subsection{Selection of Physical Fluctuations as $d\to 0$}

To have a more complete picture, it is also interesting to ask what is the limit of 
solutions of (\ref{1d-curvature}) as $d\to 0^{+}$ (the vanishing curvature limit). 
When $d=0$, equation (\ref{ge1b}) becomes the inviscid G-equation
$$
G_{t}^{\epsilon} + V({x\over \epsilon})\cdot DG^{\epsilon} + |DG^{\epsilon}|=0.
$$
It is proved in \cite{XY2010} and \cite{CNS} independently that there exists a unique $\overline H_0(p)$ such that the corresponding cell problem 
\be\label{invcell}
|p+Dw|+V(y)\cdot (p+Dw)=\overline H_0(p)  \quad \text{in $\Rset^n$}
\ee
admits a periodic (approximate) viscosity solution. This implies that 
$$
\lim_{\epsilon \to 0} G^{\epsilon}(x,t)=\bar G(x,t) \quad \text{locally uniformly in $\Rset\times [0, +\infty)$}.
$$
As in the curvature case, here $\bar G$ is the unique solution to the following effective equation, which captures the propagation of the mean flame front:
$$
\begin{cases}
\bar G_{t}+\overline H_0(D\bar G)=0\\
\bar G(x,0)=G_0(x)  \quad \text{initial flame front}.
\end{cases}
$$
The formal two-scale expansion says that 
$$
G_{\epsilon}(x,t)=\bar G(x,t)+\epsilon w(x, {x\over \epsilon})+O(\epsilon^2),
$$
where the fluctuation $w(x, \cdot)$ is a solution to (\ref{invcell}) with $p=D\bar G(x,t)$ for fixed $(x,t)$. Nevertheless, solutions to (\ref{invcell}) are in general not unique even up to a constant. This motivates

\medskip

\nit {\bf Q3:} {\it which solution to (\ref{invcell}) is the physical solution that captures the fluctuation of flame front?}

\medskip

One natural approach is to look at the limit of solutions to (\ref{curcell}) (if it exists 
uniquely) as $d\to 0$. The limit is however, 
very challenging and unknown in general. In this paper, we identify 
the limit for the equation (\ref{1d-curvature}) under some non-degeneracy conditions. 

It is easy to show that as $d\to 0^{+}$, the solution $w$ to (\ref{1d-curvature}), 
up to a subsequence, converges to a periodic viscosity solution $w_0$ of 
\be\label{cell-inviscid-1d}
\sqrt {\gamma^2+(\mu+w_{0}^{'})^2}+\gamma v(y)=\overline H_0(p) \quad \text{in $\Rset$}.
\ee
When $\gamma=0$,  $w=w_0\equiv 0$.  Without loss of generality, we set $\gamma=1$ in this section and denote $$
\overline H_0(\mu)=\overline H_0(p).
$$
Without loss of generality, in this section,  we also assume that 
$$
\max_{\Rset}v=0.
$$

\medskip

\subsubsection{Uniqueness Case}
If $|\mu|\geq \int_{0}^{1}\sqrt{(1-v)^2-1}\,dy$, $\overline H(\mu)\geq 1$ is the unique number such that 
$$
|\mu|= \int_{0}^{1}\sqrt {(\overline H(\mu)-v(y))^2-1}\,dy.
$$
Also, the inviscid equation (\ref{cell-inviscid-1d}) has a unique solution up to a constant, i.e., 
$$
w_0(x)=(sign (\mu))\int_{0}^{x}\sqrt {(\overline H(\mu)-v(y))^2-1}\,dy-\mu x+c 
$$
for some $c\in \Rset$ since $w_{0}^{\prime}+\mu$ can not change signs. Accordingly,  by $w(0)=0$, 
$$
\lim_{d\to 0^{+}}w=(sign (\mu))\int_{0}^{x}\sqrt {(\overline H(\mu)-v(y))^2-1}\,dy-\mu x.
$$

\subsubsection{Non-uniqueness Case} When $|\mu|<\int_{0}^{1}\sqrt{(1-v)^2-1}\,dy$,  $\overline H_d(\mu)=1$. The limiting problem is more interesting since solutions to the inviscid equation (\ref{cell-inviscid-1d}) are not unique if the set
$$
\mathcal{M}_0=\{x\in [0,1)|\ v(x)=\max_{\Rset}v=0\}
$$
has multiple points. For example, assume that $x_i\in \mathcal{M}_0$ for $i=1,2$. Choose $x_{\mu,i}\in (x_i,x_i+1)$ such that 
$$
\int_{x_i}^{x_{\mu,i}}\sqrt{(1-v)^2-1}\,dy-\int_{x_{\mu,i}}^{x_i+1}\sqrt{(1-v)^2-1}\,dy=\mu.
$$
Then 
$$
w_{i}(x)=
\begin{cases}
\int_{x_i}^{x}\sqrt {(1-v(y))^2-1}\,dy-\mu x,\;  \forall x\in [x_i,x_{\mu,i}]\\
\int_{x_i}^{x_{\mu,i}}\sqrt {(1-v(y))^2-1}\,dy-\int_{x_{\mu,i}}^{x}\sqrt {(1-v(y))^2-1}\,dy-\mu x,  \\
\; \forall x\in [x_{\mu,i},x_i+1]
\end{cases}
$$
(extended periodically) are both viscosity solutions to (\ref{cell-inviscid-1d}) and $w_1-w_2$ is not a constant. So a very interesting problem is to identify the solution selected by 
the limiting process, i.e., the physical fluctuation associated with the 
inviscid G-equation model. Hereafter, we assume that 
\be\label{distinct}
\text{ $\mathcal {M}_0$ is finite and $v''(x)$ is distinct for $x\in \mathcal{M}_0$}   
\ee
Choose the unique $\bar x\in \mathcal{M}_0$ such that 
$$
v''(\bar x)=\min_{x\in \mathcal{M}_0}\{-v''(x)\}.
$$
Choose $x_{\mu}\in (\bar x, \bar x +1)$ such that 
$$
\int_{\bar x}^{x_{\mu}}\sqrt{(1-v)^2-1}\,dy-\int_{x_{\mu}}^{\bar x+1}\sqrt{(1-v)^2-1}\,dy=\mu.
$$
Clearly, such $x_{\mu}$ is unique. The following is our selection result.

\begin{theo}\label{main2}
$$
\lim_{d\to 0^{+}}w=w_0(x)-w_0(0) \quad \text{uniformly in $\Rset$}.
$$
Here
\be\label{physical-sol}
w_0(x)=
\begin{cases}
\int_{\bar x}^{x}\sqrt{(1-v)^2-1}\,dy-\mu x, \; \forall x\in [\bar x, x_{\mu}]\\
\int_{\bar x}^{x_{\mu}}\sqrt{(1-v)^2-1}\,dy-\int_{x_{\mu}}^{x}\sqrt{(1-v)^2-1}\,dy-\mu x, \\
\forall x\in [x_{\mu}, \bar x+1].
\end{cases}
\ee
\end{theo}

We would like to point out that selection problems of similar spirit 
have been studied for the vanishing viscosity limit 
(\cite{LKM}, \cite{A2004}, \cite{AIPS2005}, etc), after which the viscosity solution was originally named. In these references, the 
authors aim to identify $\lim_{\epsilon\to 0^+}v_{\epsilon}$. Here $v_{\epsilon}$ is 
the unique smooth solution to 
$$
-\epsilon\Delta v_{\epsilon}+H(p+Dv_{\epsilon},x)=\overline H(p,\epsilon) \quad \text{in $\Rset^n$}.
$$
The most important case is the mechanical Hamiltonian $H(p,x)=|p|^2+G(x)$ 
with a potential function $G$. The limiting process resembles the passage  
from quantum mechanics to classical mechanics (\cite{A2004}, \cite{Evans}). 
The works \cite{A2004} and \cite{AIPS2005} deal with some special cases in high dimensions 
by employing advanced tools from dynamical systems and random perturbations. 
Assumptions therein are very hard to check however. 
The method in \cite{LKM} is purely 1-d. Based on simple comparison principles of PDEs/ODEs, our arguments are simpler and more robust. In particular, they can be easily extended to handle certain cases in high dimensions. 
The rest of the paper contains the proofs of the main theorems.

\section{Proof of Theorem \ref{main1}}

\nit Proof: (1) is trivial. Let us prove (2) which is the most difficult and interesting part. Fix $(\gamma, \mu)$. Denote $\phi={\mu+w'\over \gamma}$.Then $\phi$ is the unique periodic solution to
$$
-{d\phi'\over 1+\phi^2}+\sqrt {1+\phi^2}+v(y)=E(d)={\overline H_{d}(p)\over \gamma}  \quad \text{in $\Rset$}
$$
 subject to $\int_{0}^{1}\phi(x)\,dx={\mu\over \gamma}$ . To prove (2) is equivalent to showing that 
$$
E'(d)<0.
$$
Taking derivative on both sides of the above equation with respect to $d$, we obtain that 
$$
-dF'+b(x)F=E'(d)(1+\phi^2)+\phi',
$$
where  $b(x)={2d\phi'\phi\over 1+\phi^2}+\phi \sqrt {1+\phi^2}$ and $F(x)=\phi_d(x)$, i.e., the derivative of $\phi$ with respect to $d$. Clearly, $F$ is periodic and has zero mean, i.e., $\int_{[0,1]}F=0$. Note that $v$ is not constant is equivalent to saying the $\phi$ is not constant. Then (2) follows immediately from Lemma \ref{key-lemma}. 

\medskip

\nit (3) Integrating both sides of (\ref{1d-curvature}), we obtain:
\be\label{average}
\overline H_d(p)=\int_{0}^{1}\sqrt {\gamma^2+(\mu+w')^2}\,dy+\gamma \int_{0}^{1}v(y)\,dy.
\ee
So due to the convexity of $s(t)=\sqrt{\gamma^2+t^2}$, 
$$
\overline H_d(p)\geq |p|+\gamma \int_{0}^{1}v(y)\,dy.
$$
Also, by maximum principle, we have that 
$$
\overline H_d(p)\leq |p|+\max_{\Rset} \gamma v
$$
and
$$
\max_{\Rset}|\mu+w'|\leq \overline H_d(p)-\min_{\Rset} \gamma v\leq |p|+2\max_{\Rset} |\gamma v|.
$$
Hence, up to a sequence, we may assume that 
$$
\lim_{d\to 0}\overline H_d=\overline H_0 \quad \mathrm{and} \lim_{d\to 0^{+}}w=w_0 \quad \text{uniformly in $\Rset$}.
$$
Then the stability of viscosity solution immediately implies that $w_0$ is a continuous periodic viscosity solution to 
$$
\sqrt {\gamma^2+\left(\mu+w_{0}^{'}\right)^2}+\gamma v(y)=\overline H_0(p) \quad \text{in $\Rset$}. 
$$
Note that $\overline H_0(p) $ is unique number such that the above equation has a periodic viscosity solutions $w_0$ although $w_0$ might not be unique. See \cite{LPV} for general cases.

\medskip

\nit (4). If $\gamma=0$, this is trivial. So we assume that $\gamma\not=0$. Note that estimates of $\overline H_d$ and $\mu+w'$ in (3) are independent of $d$. Since 
$$
w''={1\over d\gamma^2}(\gamma^2+(\mu+w')^2)\left(\sqrt{\gamma^2+(\mu+w')^2}+v-\overline H_d(\mu)\right),
$$
we have that 
$$
\max_{\Rset}|w''|\leq {C\over {d}}
$$
for a constant $C$ independent of $d$.  Due to the periodicity of $w$ and $w(0)=0$,   it is obvious that
$$
\lim_{d\to +\infty}w=\lim_{d\to +\infty}w'=0 \quad \text{uniformly in $\Rset$}.
$$
Combining with (\ref{average}), (4) holds.

\qed

\begin{lem}\label{key-lemma} Let $d>0$ and $\phi$ be a non-constant $C^1$ periodic function. If the following equation has a mean-zero, periodic solution  $F$
$$
-dF'+b(x)F=\phi'+\alpha (1+\phi^2)  \quad \text{in $\Rset$}
$$
for some $\alpha\in \Rset$ and 
$$
b(x)={2d\phi'\phi\over 1+\phi^2}+\phi \sqrt {1+\phi^2}, 
$$
then
$$
\alpha<0.
$$
\end{lem}
Proof:  It suffices to prove this for $d=1$. The proof for other $d$ 
is similar.  We can solve $F$ in terms of $\phi$ and $\alpha$. Using $F$ is periodic and mean zero (i.e., $F(0)=F(1)$ and $\int_{0}^{1}F(s)\,ds=0$), it is easy to obtain that 
$$
\alpha=-{e^{g(1)}\int_{0}^{1}\phi'e^{-g(x)}\,dx\int_{0}^{1}e^{g(x)}\,dx-(e^{g(1)}-1)\int_{0}^{1}e^{g(x)}\int_{0}^{x}\phi'e^{-g(y)}\,dydx \over e^{g(1)}\int_{0}^{1} (1+\phi^2)e^{-g(x)}\,dx\int_{0}^{1}e^{g(x)}\,dx-(e^{g(1)}-1)\int_{0}^{1}e^{g(x)}\int_{0}^{x} (1+\phi^2)e^{-g(y)}\,dydx}.
$$
Here 
$$
g(x)=\int_{0}^{x}b(y)\,dy=\log (1+\phi^2(x))-\log (1+\phi^2(0))+\int_{0}^{x}\phi\sqrt{1+\phi^2}\,dx.
$$
In particular, $g(1)=\int_{0}^{1}\phi\sqrt{1+\phi^2}\,dx$. The denominator is obviously positive. Hence $\alpha<0$ is equivalent to proving the inequality
$$
e^{g(1)}\int_{0}^{1}\phi'e^{-g(x)}\,dx\int_{0}^{1}e^{g(x)}\,dx>(e^{g(1)}-1)\int_{0}^{1}e^{g(x)}\int_{0}^{x}\phi'e^{-g(y)}\,dydx 
$$
for every non-constant $C^1$ periodic function $\phi$. Denote that 
$$
h(x)=\int_{0}^{x}\phi\sqrt{1+\phi^2}\,dy.
$$
Then it is equivalent to showing that 
$$
e^{h(1)}\int_{0}^{1}{\phi'\over 1+\phi^2}e^{-h(x)}\,dx\int_{0}^{1}(1+\phi^2)e^{h(x)}\,dx>(e^{h(1)}-1)\int_{0}^{1}(1+\phi^2)e^{h(x)}\int_{0}^{x}{\phi'\over 1+\phi^2}e^{-h(y)}\,dy
$$
Write $\lambda(\phi)=\arctan \phi$. Using integration by parts and $\phi(0)=\phi(1)$, we have that 
$$
LHS=e^{h(1)}\left(\lambda(\phi(1))e^{-h(1)}-\lambda(\phi(1))+\int_{0}^{1}\lambda(\phi)e^{-h(x)}\phi\sqrt {1+\phi^2}\,dx\right)\int_{0}^{1}(1+\phi^2)e^{h(x)}\,dx.
$$
and the RHS is 
$$
\begin{array}{ll}
RHS=&(e^{h(1)}-1)\left(\int_{0}^{1}\lambda(\phi)(1+\phi^2)\,dx-\lambda(\phi(1))\int_{0}^{1}(1+\phi^2)e^{h(x)}\,dx\right)\\[3mm]
&+(e^{h(1)}-1)\left(\int_{0}^{1}(1+\phi^2)e^{h(x)}\int_{0}^{x}\lambda(\phi)e^{-h(y)}\phi\sqrt{1+\phi^2}\,dydx\right).
\end{array}
$$
By Fubini Theorem, 
\ba
& & \int_{0}^{1}(1+\phi^2)e^{h(x)}\int_{0}^{x}\lambda(\phi)e^{-h(y)}\phi\sqrt{1+\phi^2}\,dydx \no \\
& & =\int_{0}^{1}\lambda(\phi)e^{-h(x)}\phi\sqrt{1+\phi^2}\int_{x}^{1}(1+\phi^2)e^{h(y)}\,dydx. \no 
\ea

Then $LHS-RHS$ is $A+B-C$ for 
$$
A(\phi)=e^{h(1)}\int_{0}^{1}\lambda(\phi)e^{-h(x)}\phi\sqrt {1+\phi^2}\int_{0}^{x}(1+\phi^2)e^{h(y)}\,dydx,
$$
$$
B(\phi)=\int_{0}^{1}\lambda(\phi)e^{-h(x)}\phi\sqrt{1+\phi^2}\int_{x}^{1}(1+\phi^2)e^{h(y)}\,dydx.
$$
and
$$
C(\phi)=(e^{h(1)}-1)\int_{0}^{1}\lambda(\phi)(1+\phi^2)\,dx. 
$$
If $h(1)=0$, then $A+B-C=A+B\geq 0$ since $s\lambda (s)\geq 0$. Cleary, $ ``=0"$ if and only if $\phi\equiv 0$.  So we assume that 
$$
h(1)\not=0.
$$
Also, note that for $\tilde \phi(x)=-\phi(-x)$, the correspsonding
$$
\tilde b(x)={2\tilde \phi'\tilde \phi\over 1+{\tilde \phi}^2}+\tilde \phi \sqrt {1+{\tilde \phi}^2}=-b(-x)
$$
and $\tilde F(x)=-F(-x)$ satisfies that 
$$
-\tilde F'+\tilde b(x)\tilde F=\tilde \phi'+\alpha(1+{\tilde \phi}^2).
$$
Hence, without lost of generality,  we may further assume that 
$$
h(1)>0.
$$
Denote $\phi_{+}=\max\{\phi, 0\}$ and $\phi_{-}=\min\{\phi,0\}$. Aso write
$$
h^{\pm}(x)=\int_{0}^{x}\phi_{\pm}\sqrt{1+\phi_{\pm}^{2}}\,dy. 
$$
Note that $h(x)=h^{+}+h^{-}$. Now let us prove the following lemma. 
\begin{lem}\label{sign} We have that 
$$
A(\phi)+B(\phi)-C(\phi)\geq e^{h^{-}(1)}\left(A(\phi_{+})+B(\phi_{+})-C(\phi_{+})\right).
$$
The equality holds if only if $\phi\geq 0$, i.e., $\phi_{-}=0$.
\end{lem}

\nit Proof:  Clearly
$$
\begin{array}{ll}
A(\phi)&\geq e^{h(1)}\int_{0}^{1}\lambda(\phi_{+})e^{-h(x)}\phi_{+}\sqrt{1+\phi_{+}^{2}}\int_{0}^{x}(1+\phi_{+}^{2})e^{h(y)}\,dydx\\[5mm]
&=e^{h(1)}\int_{0}^{1}\lambda(\phi_{+})e^{-h^{+}(x)}\phi_{+}\sqrt{1+\phi_{+}^{2}}\int_{0}^{x}(1+\phi_{+}^{2})e^{h^{+}(y)}e^{h^{-}(y)-h^{-}(x)}\,dydx\\[5mm]
&\geq e^{h^{-}(1)}A(\phi_{+}), \quad \text{since $h^{-} (x)\leq h^{-} (y)$ for $x\geq y$}.
\end{array}
$$
Also, 
$$
\begin{array}{ll}
& B(\phi)\geq \int_{0}^{1}\lambda(\phi_{+})e^{-h(x)}\phi_{+}\sqrt{1+\phi_{+}^{2}}\int_{x}^{1}(1+\phi_{+}^{2})e^{h(y)}\,dydx\\[5mm]
&=e^{h^{-}(1)}\int_{0}^{1}\lambda(\phi_{+})e^{-h^{+}(x)}\phi_{+}\sqrt{1+\phi_{+}^{2}}\int_{x}^{1}(1+\phi_{+}^{2})e^{h^{+}(y)}e^{h^{-}(y)-h^{-}(1)}e^{-h^{-}(x)}\,dydx\\[5mm]
&\geq e^{h^{-}(1)}B(\phi_{+}) \quad \text{since $0\geq h^{-}(y)\geq h^{-}(1)$ for all $y\in [0,1]$}
\end{array}
$$
and
$$
\begin{array}{ll}
C(\phi)&\leq (e^{h(1)}-1)\int_{0}^{1}\lambda(\phi_{+})(1+\phi_{+}^{2})\,dx\\[5mm]
&={ (e^{h(1)}-1)\over (e^{h^{+}(1)}-1)} C(\phi_{+})\\[5mm]
&\leq e^{h^{-}(1)}C(\phi_{+}).
\end{array}
$$
Obviously, for all inequalities to hold, we must have $h^{-}\equiv 0$ and $\phi_{-}\equiv 0$. \qed

Now let us  continue the proof of Lemma  \ref{key-lemma}.  Since $h(1)>0$, that $\phi$ is not constant implies $\phi_+$ is not constant either. By a small perturbation  like $\phi_{+}+\epsilon$, we may assume that $\phi_{+}>0$ in computations below. Then $h^{+}$ is strictly increasing.  After changing of variables $h^{+}(x)\to x$ and writing $\psi(h^{+}(x))=\phi_{+}(x)$ and $T=h^{+}(1)$, we obtain that 
$$
A(\phi_{+})=A_{T, \psi}=e^{T}\int_{0}^{T}\lambda(\psi)e^{-x}\int_{0}^{x}{\sqrt {1+{\psi^2}}\over \psi }e^y\,dydx,
$$
$$
B(\phi_{+})=B_{T, \psi}=\int_{0}^{T}\lambda(\psi)e^{-x}\int_{x}^{T}{\sqrt {1+{\psi^2}}\over \psi }e^y\,dydx
$$
and
$$
C(\phi_{+})=C_{T, \psi}=(e^T-1)\int_{0}^{T}\lambda(\psi){\sqrt {1+{\psi^2}}\over \psi }\,dx.
$$
So 
\ba
A_{T,\psi}+B_{T,\psi}-C_{T,\psi} &=& 
e^{T}\int_{0}^{T}\lambda(\psi)e^{-x}\int_{0}^{x}{\sqrt {1+{\psi^2}}\over \psi }\, e^y\,dydx
\no \\
&+ & \int_{0}^{T}\lambda(\psi)e^{-x}\int_{x}^{T}{\sqrt {1+{\psi^2}}\over \psi }e^y\,dydx \no \\
&-& (e^T-1)\int_{0}^{T}\lambda(\psi){\sqrt {1+{\psi^2}}\over \psi }\,dx. \no 
\ea
Let $M=\max_{[0,T]}\psi=\max_{[0,1]}\phi_{+}>0$. According to Theorem \ref{continuous} by taking $f(x)=\lambda(\psi)=\arctan(\psi)$, $g(y)={1\over \sin y}$, $L=\arctan(M)$ and $\theta={1\over \sqrt {1+M^2}}$, we  have that ${\sqrt {1+{\psi^2}}\over \psi }=g(f)$ and 
$$
\begin{array}{ll}
A_{T,\psi}+B_{T,\psi}-C_{T,\psi}&\geq {1\over 2\sqrt {1+M^2}}\int_{[0,T]^2}|\lambda(\psi (x))-\lambda(\psi (y))|^2\,dxdy\\[5mm]
&={1\over 2\sqrt {1+M^2}}\int_{[0,1]^2}|\lambda(\phi_{+} (x))-\lambda(\phi_{+} (y))|^2J(x)J(y)\,dxdy\\[5mm]
&>0 \quad \text{since $\phi_+$ is not constant}.
\end{array}
$$
Here $J(x)=\phi_{+}(x)\sqrt{1+\phi_{+}^{2}}$. Combining with  Lemma \ref{sign}, $A(\phi)+B(\phi)-C(\phi)>0$. \qed

\subsection{The Key Inequalities}

Given $n\in \Nset$. Let $\{b_{ik}\}_{1\leq i,k\leq n}$ and $\{\tilde b_{ik}\}_{1\leq i, k\leq n}$ be two given sequences of positive numbers satisfying that for all $i,k$
$$
\sum_{l=1}^{i}b_{il}+\sum_{l=i}^{n}\tilde b_{il}=\sum_{l=k}^{n}b_{lk}+\sum_{l=1}^{k}\tilde b_{lk}=c.
$$
Here $c$ is a constant independent of $i$ and $k$.  Also, 
\be\label{lowbound}
\min\{\min_{1\leq k\leq i\leq n} b_{ik}, \min_{1\leq i\leq k\leq n}\tilde b_{ik}\}\geq \tau>0.
\ee

\begin{lem}\label{discrete} Assume that $L>0$ and $g\in C((0,L])$ satisfies

$$
g'(a)\leq -\theta  \quad \text{ for some $\theta\geq 0$}.
$$
\nit Then
$$
\sum_{i=1}^{n}a_i\sum_{k=1}^{i}g(a_k)b_{ik}+\sum_{i=1}^{n}a_i\sum_{k=i}^{n}g(a_k)\tilde b_{ik}\geq c\sum_{i=1}^{n}a_ig(a_i)+{\theta\tau\over 2}\sum_{1\leq i,k\leq n}(a_i-a_k)^2.
$$
for all $(a_1,a_2,...,a_n)\in (0,L]^n$.  Here $\tau$ is from (\ref{lowbound}). Moreover, if $\theta>0$, the equality holds if and only if $a_1=a_2=..=a_n$. 
\end{lem}

\nit Proof:  By approximation, we may assume that $\theta>0$. For convenience, denote
$$
W(a_1,a_2,...,a_n)=\sum_{i=1}^{n}a_i\sum_{k=1}^{i}g(a_k)b_{ik}+\sum_{i=1}^{n}a_i\sum_{k=i}^{n}g(a_k)\tilde b_{ik}
$$
and 
$$
H(a_1,a_2,...,a_n)=c \sum_{i=1}^{n}a_ig(a_i)+{\theta\tau\over 2}\sum_{1\leq i,k\leq n}(a_i-a_k)^2.
$$
It suffices to show that for any fixed $r\in (0,L)$, 
$$
\min_{[r,L]^n}(W-H)=0
$$
and the minimum is attained when all $a_i$ are the same. 

Choose  $(\hat a_1,\hat a_2,\hat a_3,..\hat a_n)\in [r,L]^n$ such that 
$$
W(\hat a_1,\hat a_2,\hat a_3,..\hat a_n)-H(\hat a_1,\hat a_2,\hat a_3,..\hat a_n)=\min_{[r,L]^n}(W-H).
$$
Assume that $\hat a_j=\max_{1\leq i\leq n}\{\hat a_i\}$.  If $\hat a_j=r$, then $\hat a_1=\hat a_2=..=\hat a_n=r$ and we are done.  So let us assume that 
$$
\hat a_j>r.
$$
Then
$$
W_{a_j}-H_{a_j}\leq 0 \quad \text{at $(\hat a_1,\hat a_2,\hat a_3,..\hat a_n)$}.
$$
Here we include $<0$ since $\hat a_j$ might be equal to $L$. Accordingly, 
$$
\begin{array}{ll}
&\sum_{k=1}^{j}g(\hat a_k)b_{jk}+\sum_{k=j}^{n}g(\hat a_k)\tilde b_{jk}\\[5mm]
&+g'(\hat a_j)\sum_{k=j}^{n}\hat a_kb_{kj}+g'(\hat a_j)\sum_{k=1}^{j}\hat a_k\tilde b_{kj}\\[5mm]
&\leq c (g(\hat a_j)+\hat a_jg'(\hat a_j))+2\sum_{k\not= j}\theta \tau (\hat a_j-\hat a_k)
\end{array}
$$
On the other hand, since  $g'\leq -\theta<0$, we also have that 
$$
\begin{array}{ll}
&\sum_{k=1}^{j}g(\hat a_k)b_{jk}+\sum_{k=j}^{n}g(\hat a_k)\tilde b_{jk}\\[5mm]
&+g'(\hat a_j)\sum_{k=j}^{n}\hat a_kb_{kj}+g'(\hat a_j)\sum_{k=1}^{j}\hat a_k\tilde b_{kj}\\[5mm]
&\geq c(g(\hat a_j)+\hat a_jg'(\hat a_j))+2\sum_{k\not= j}\theta \tau (\hat a_j-\hat a_k).
\end{array}
$$
Hence all equalities should hold and $\hat a_1=\hat a_2...=\hat a_n$ follows from that $g$ is strictly decreasing. Then $W(\hat a_1,\hat a_2,\hat a_3,..\hat a_n)-H(\hat a_1,\hat a_2,\hat a_3,..\hat a_n)=0$.\qed

Now we are ready to state a specific continuous version for our purpose.

\begin{theo}\label{continuous} Let $T>0$ and $f\in C([0,T])$  be  a continuous positive function. Suppose that $g\in C^1((0,L])$ for $L=\max_{[0,T]}f$. 

\nit (1) If  $g'\leq -\theta$ for some $\theta\geq 0$, then
$$
\begin{array}{ll}
&e^{T}\int_{0}^{T}f(x)e^{-x}\int_{0}^{x}g(f(y))e^y\,dydx+\int_{0}^{T}f(x)e^{-x}\int_{x}^{T}g(f(y))e^y\,dydx\\[5mm]
&\geq (e^T-1)\int_{0}^{T}f(x)g(f(x)))\,dx+{\theta\over 2} \int_{[0,T]^2}|f(x)-f(y)|^2\,dxdy.
\end{array}
$$
(2) If If  $g'\geq \theta$ for some $\theta\geq 0$, then

$$
\begin{array}{ll}
&e^{T}\int_{0}^{T}f(x)e^{-x}\int_{0}^{x}g(f(y))e^y\,dydx+\int_{0}^{T}f(x)e^{-x}\int_{x}^{T}g(f(y))e^y\,dydx\\[5mm]
&\leq (e^T-1)\int_{0}^{T}f(x)g(f(x)))\,dx-{\theta\over 2} \int_{[0,T]^2}|f(x)-f(y)|^2\,dxdy.
\end{array}
$$

\end{theo}

\nit Proof.  (1) For $n\in N$, let $x_i={iT\over n}$ for $i=1,2,..,n$. Note that for $i,k=1,2,3,..n$, 
$$
\sum_{l=1}^{i}e^{T-x_i+x_l}+\sum_{l=i}^{n}e^{x_l-x_i}={e^{T+{T\over n}}-1\over e^{T\over n}-1}=\sum_{l=k}^{n}e^{T-x_l+x_k}+\sum_{l=1}^{k}e^{x_k-x_l}.
$$
Then desired inequality in (1)  follows from Lemma \ref{discrete} and Riemann sum approximation by taking $a_i=f(x_i)$, $c={e^{T+{T\over n}}-1\over e^{T\over n}-1}$, $\tau=1$, 
$$
b_{ik}=e^{T-x_i+x_k}  \quad  \mathrm{and} \quad \tilde b_{ik}=e^{x_k-x_i} \quad \text{for $1\leq i, k\leq n$}. 
$$

(2) follows immediately from (1) by considering $-g$. 
\qed

\begin{rmk}\label{vanishing-viscosity} Similar to the proof of Theorem \ref{main1}, (1) in the above  Theorem \ref{continuous} also implies that the one dimensional viscous effective Hamiltonian $\overline H(p,d)$ given by the cell problem
$$
-d\, w''+H(p+w')+G(x)=\overline H(p,d) \quad \text{in $\Rset$}
$$
is strictly decreasing with respect to the diffusivity $d>0$ for a non-constant function $G$, and 
a strictly convex function $H:\Rset\to \Rset$. Here we choose $f=p+w'$ and $g={1\over H'}$ after suitable translations.  It remains an interesting problem whether this is also true in high dimensions. For the special case $H(p)={1\over 2}|p|^2$,  using integration by parts, it is easy to derive that 
\be\label{viscous-basic}
{\partial \overline H(p,d)\over \partial d}=-{\int_{\Bbb T^n}|Dw_d|^2e^{-w_d}\,dx\over \int_{\Bbb T^n}e^{-w_d}\,dx}\leq 0
\ee
and $``="$ holds if and only if $G$ is a constant. Here  $w_d$ represents the derivative of $w$ with respect to $d$. On the other hand,  if $H$ is non-convex, then (2) in the above  Theorem \ref{continuous} implies that for some $p$,  $\overline H(p,d)$ could be strictly increasing with respect to $d$. 
\end{rmk}

\subsection{Calculations in High Dimensions in Perturbative Cases.}

Consider the case of weak flow or $\delta V$ for $0\leq \delta\ll 1$. Let $p\in \Rset^n$ be a unit vector satisfying the Diophantine condition, i.e., there exist $\beta, C>0$ such that
$$
|p\cdot \vec{k}|\geq {C\over |\vec{k}|^{\beta}} \quad \text{for all $\vec{k}\in \Zset^n\backslash\{0\}$}.
$$
Owing to \cite{LS2005}, when $\delta$ is small eough,   the cell problem  (\ref{curcell}) has a viscosity solution. Formally, we can write the solution as 
$$
w=\delta w_1+\delta^2 w_2+O(\delta^3)
$$
and the constant (turbulent flame speed)
\be\label{barH-expansion}
\overline H_d(p)=|p|+\delta \alpha_1(p)+\delta^2 \alpha_2(p)+O(\delta^3).
\ee
By comparing coefficients of $\delta$ and $\delta^2$, $w_1$ and $w_2$ are 
determined by inhomogeneous linear equations. They can be solved in terms of Fourier series. 
For example, $w_1$ satisfies 
$$
-d(\Delta w_1-p\cdot D^2w_1\cdot p)+p\cdot Dw_1+p\cdot V=\alpha_1(p).
$$
The equation for $w_2$ is more messy. Applying Fredholm alternatives to both equations, we have that 
$$
\alpha_1(p)=p\cdot \int_{\Bbb T^n}V\,dx=p\cdot \lambda_0
$$
and
$$
\alpha_2(p)={1\over 2}\int_{\Bbb T^n}|Dw_1|^2\,dx={1\over 2}\sum_{\vec{k}\in \Zset^n\backslash\{0\}}{|p\cdot \lambda_{\vec{k}}|^2|\vec{k}|^2\over d^2(|\vec{k}|^2-|p\cdot \vec{k}|^2)^24\pi^2+|p\cdot \vec{k}|^2},
$$
where $\lambda_{\vec{k}}\in \Cset^n$ are Fourier coefficients of $V$, i.e., $V=\sum_{\vec{k}\in \Zset^n}\lambda_{\vec{k}}e^{i2\pi \vec{k}\cdot x}$. Clearly, $\overline H_d(p)$ is strictly decreasing with respect to $d$.  The approximation of $\overline H_d(p)$ (\ref{barH-expansion}) can actually be proved easily  through maximum principles of viscosity solutions, i.e., evaluating at where $w-\delta w_1-\delta^2 w_2$ attains maximum/minimum values.

\section{Proof of Theorem \ref{main2}}

Let us first prove some lemmas. Recall that 
$$
\mathcal{M}_0=\{x\in [0,1)|\ v(x)=\max_{\Rset}v=0\}. 
$$
 See section 1.2.2 (Non-uniqueness Case) for the range of $\mu$,  defintions of $\bar x$ and $x_\mu$ and other assumptions like (\ref{distinct}).

\begin{lem}\label{slope} Assume that $\mathcal{M}_0=\{\bar x \}$, i.e., it contains a single element. Then
$$
\lim_{d\to 0^{+}}{\overline H_d(\mu)-1\over d}=-\sqrt {-v''(\bar x)}.
$$
\end{lem}

\nit Proof: Since $\mathcal{M}_0$ has only one element,  $1-v>1$ in $(\bar x, \bar x+1)$. Then it is easy to see that   periodic viscosity solutions to 
$$
\sqrt{1+(\mu+w_{0}^{'})^2}+v(y)=1  \quad \text{in $\Rset$}
$$
are unique up to a constant. Hence, since $w(0)=0$, 
\be\label{limit}
\lim_{d\to 0^{+}}w=w_0(x)-w_0(0) \quad \text{uniformly in $\Rset$}.
\ee
Here $w_0$ is given by (\ref{physical-sol}).  Fix $\delta>0$ and denote
$$
u_{\delta, \pm}(x)=
\begin{cases}
\int_{\bar x}^{x}\sqrt{(1-(1\pm \delta) v)^2-1}\,dy \quad \text{for $x\geq \bar x$}\\[3mm]
\int_{x}^{\bar x}\sqrt{(1-(1\pm \delta) v)^2-1}\,dy \quad \text{for $x\leq \bar x$}.
\end{cases}
$$
Apparently,
$$
u_{\delta, -}(x)<u_0(x)=w_0(x)+\mu x<u_{\delta, +}(x) \quad \text{for $x\in [x_{\mu}-1, x_{\mu}]\backslash \{\bar x\}$}
$$
and $u_{\delta, -}(\bar x)=u_0(\bar x)=u_{\delta, +}(\bar x)=0$. See the left picture on Figure 4.  Denote 
$$
e_{\delta}=\min_{\text{$x=x_{\mu}$ or $x_{\mu}-1$}} \{ u_0(x)-u_{\delta, -}(x),\ u_{\delta, +}(x)-u_0(x)\}>0.
$$
and
$$
u_{d,\delta, \pm}(x)=w(x)-w(\bar x)+\mu (x-\bar x)\pm {1\over 2}e_{\delta}.
$$
Clearly, by (\ref{limit}), when $d$ is small enough, there exist $x_{d,\delta, \pm}\in (x_\mu-1, x_\mu)$ such that
$$
u_{d,\delta,+}(x_{d,\delta, +})-u_{\delta,+}(x_{d, \delta, +})\geq u_{d,\delta,+}(x)-u_{\delta,+}(x) \quad \text{for all $x\in (x_\mu-1, x_\mu)$}
$$
and
$$
u_{d,\delta,-}(x_{d, \delta,-})-u_{\delta,-}(x_{d, \delta,-})\leq u_{d,\delta,-}(x)-u_{\delta,-}(x) \quad \text{for all $x\in (x_\mu-1,x_\mu)$}.
$$
Hence maximum principle implies that 
$$
-{du_{\delta,+}^{''}\over 1+(u_{\delta, +}^{'})^2}+\sqrt {1+(u_{\delta, +}^{'})^2}+v\leq \overline H_d(\mu) \quad \text{at $x_{d,\delta,+}$}.
$$
So
$$
-{du_{\delta,+}^{''}\over 1+(u_{\delta, +}^{'})^2}\leq \overline H_d(\mu)-1+\delta v\leq   \overline H_d(\mu)-1\quad \text{at $x_{d,\delta,+}$}.
$$
Sending $d\to 0$ first  and then $\delta\to 0$, we derive that  $x_{d, \delta,+}\to \bar x$ and
$$
\liminf_{d\to 0^{+}}{\overline H_d(\mu)-1\over d}\geq -\sqrt {-v''(\bar x)}.
$$
By looking at $x_{d, \delta, -}$, similarly, we can obtain that 
$$
\limsup_{d\to 0^{+}}{\overline H_d(\mu)-1\over d}\leq -\sqrt {-v''(\bar x)}.
$$
Hence we finish the proof. \qed

\begin{center}
\includegraphics[scale=0.25]{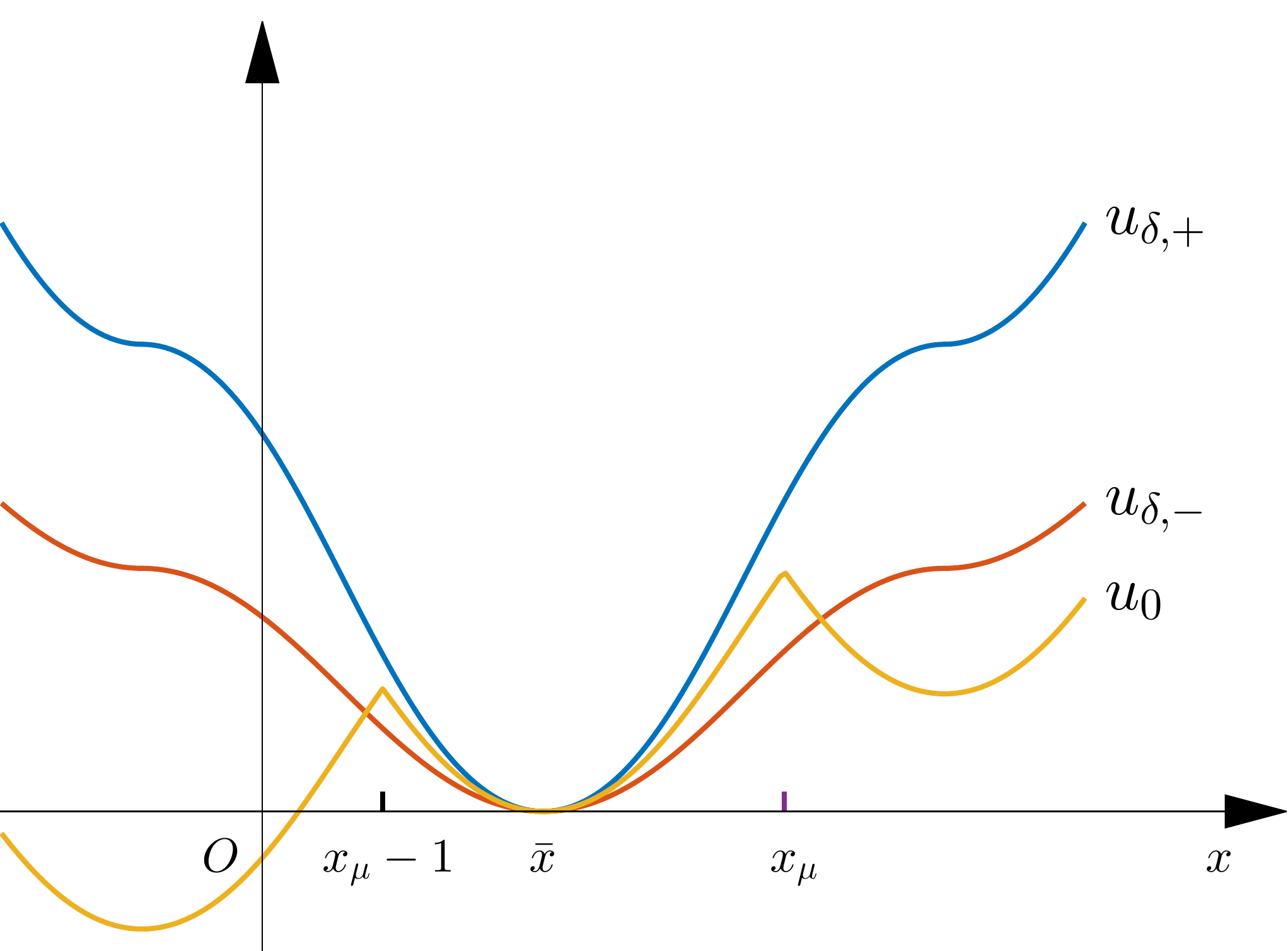} \quad \qquad
\includegraphics[scale=0.3]{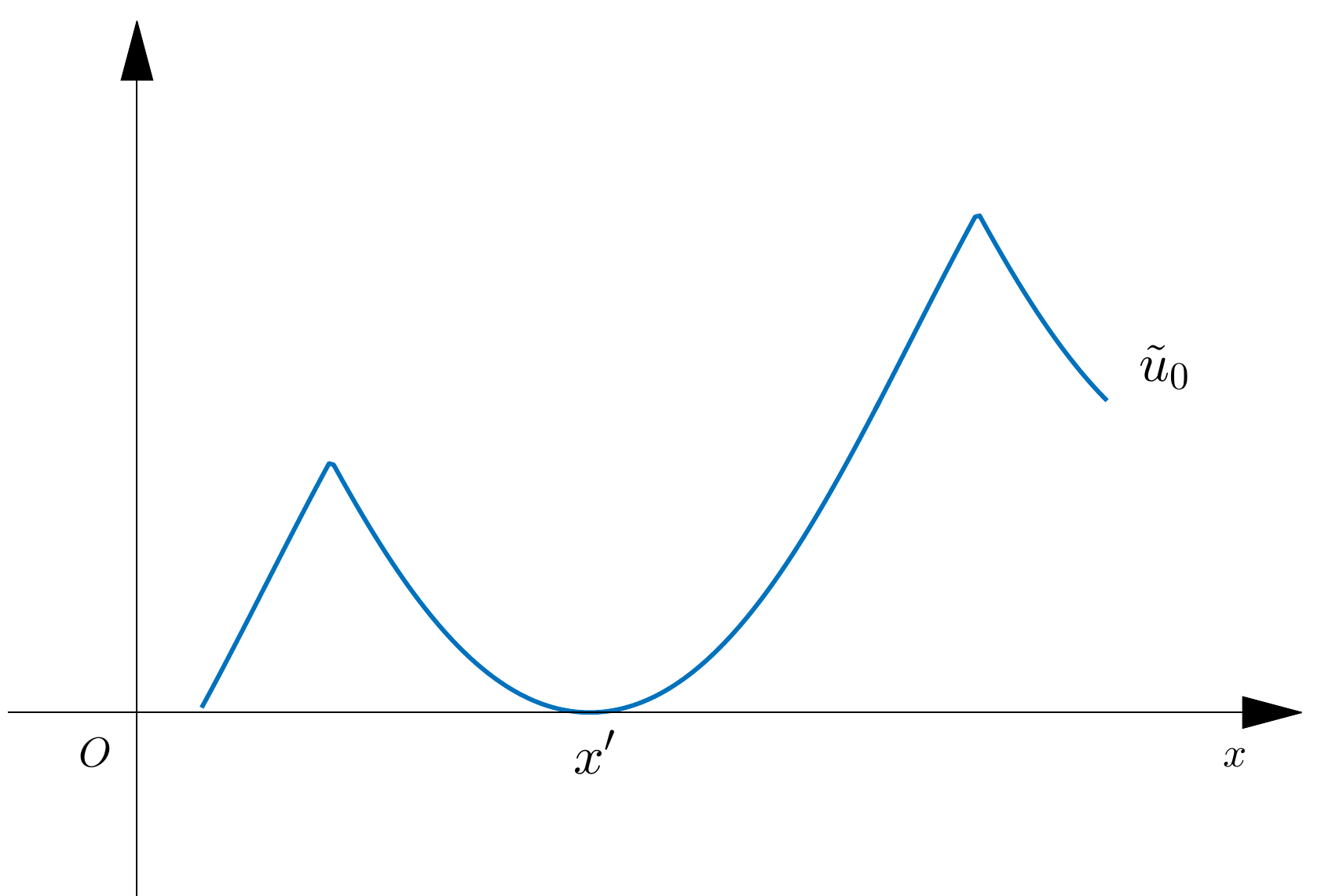}\quad \qquad
\captionof{figure}{Left: graphes of $u_{\delta,\pm}$ and $u_0$. \qquad Right: Turning points}
\end{center}

\begin{rmk}\label{u-shape}The above proof based on comparison and maximum principle actually also shows that for any subsequence $\{d_m\} \to 0$, if 
$$
\lim_{d_m\to 0^{+}}w(x)=\tilde w_{0}(x)
$$
and $\tilde u_0=\mu x+\tilde w_{0}(x)$ has turning point at some $x'\in \mathcal{M}$, i.e. there exists a $\tau>0$ such that (see the right picture on Figure 4)
$$
\tilde u_0(x)-\tilde u_0(x')=
\begin{cases}
\int_{x'}^{x}\sqrt{(1- v)^2-1}\,dy \quad \text{for $x\in [x', x'+\tau]$}\\[3mm]
\int_{x}^{x'}\sqrt{(1-v)^2-1}\,dy \quad \text{for $x\in [x'-\tau, x']$}, 
\end{cases}
$$
then
$$
\lim_{m\to +\infty}{\overline H_{d_m}(\mu)-1\over d_m}=-\sqrt{-v''(x')}.
$$
\end{rmk}

\begin{lem}\label{mini-turn} Suppose that $\tilde w$ is a periodic viscosity solution to the inviscid equation
$$
\sqrt{1+(\mu+\tilde w')^2}+v=1 \quad \text{in $\Rset$}.
$$ 
Then $x_0\in \Rset$ is a turning point of $\tilde u(x)=\mu x+\tilde w$ if and only if ${\tilde u}(x)$ attains local minimum at $x_0$. 
\end{lem}

\nit Proof: ``$\Rightarrow$" is obvious. We only need to show that any local minimum point $x_0$ must be a turning point. By the definition of viscosity solutions, 
$$
1+v(x_0)\geq 1.
$$
So $v(x_0)=0$ and $x_0\in \mathcal{M}_0$. Choose $\tau>0$ such that $(x_0, x_0+\tau)\cap \mathcal{M}_0=\emptyset$ and ${\tilde u}'(x_0+\tau)=p+\tilde w'(x_0+\tau)>0$. Then we must have that 
$$
{\tilde u}'(y)>0 \quad \text{for any $y\in (x_0,x_0+\tau)$ where ${\tilde u}'$ exists.}
$$
Otherwise there will be a local mimimum point in $(x_0,x_0+\tau)$. Note that any local minimum point belongs to $\mathcal{M}_0$. This will contradict to the choice of $\tau$. Accordingly, 
$$
{\tilde u}'=\sqrt {(1-v)^2-1} \quad \text{in $(x_0,x_0+\tau)$}.
$$
Similarly, we can show that for some $\tau'>0$,
$$
{\tilde u}'=-\sqrt {(1-v)^2-1} \quad \text{in $(x_0-\tau',x_0)$.}
$$
\qed

\medskip

\nit {\bf Proof of Theorem \ref{main2}}.

\nit {\bf Step 1:} We first show that 
\be\label{one-limit}
\liminf_{d\to 0^{+}}{\overline H_d(\mu)-1\over d}\geq -\sqrt {-v''(\bar x)}.
\ee

In fact, let $h(x)$ be a smooth periodic function satisfying that $h(\bar x)=0$ and $h(x)>0$ for $x\notin \bar x+\Zset$. For $\epsilon>0$, denote 
$$
v_{\epsilon}(x)=v(x)-\epsilon h(x)
$$
and $\overline H_{d,\epsilon}(p)$ from the cell problem (\ref{1d-curvature}) with $\gamma=1$ and $v$ replaced by $v_{\epsilon}$. It is easy to see that 
$$
\overline H_{d}(\mu)\geq \overline H_{d,\epsilon}(\mu).
$$
Choose $\epsilon$ small enough such that
$$
|\mu|<\int_{0}^{1}\sqrt{(1-v_{\epsilon})^2-1}\,dx.
$$
Clearly,  $\max_{\Rset}v_{\epsilon}=0$ and the maximum is only obtained at $\bar x+\Zset$. Then (\ref{one-limit}) followes immediately from Lemma \ref{slope}.

\medskip

\nit {\bf Step 2:} Suppose $\tilde u=\mu x+\tilde w$ is the limit of a subsequence of $\mu x+w$ as $d\to 0$. Combining with the above Remark \ref{u-shape} and assumption (\ref{distinct}), (\ref{one-limit}) implies that  $\tilde u$  can only have a turning point at $\bar x$. Owing to Lemma \ref{mini-turn}, $\tilde u$ does not have local minimum points in $(\bar x, \bar x+1)$.  Together with
$|\mu|<\int_{0}^{1}\sqrt{(1-v)^2-1}\,dx$,  it is easy to see that there exists a unique $x_\mu\in (\bar x, \bar x+1)$ such that $\tilde u$ is increasing in $(\bar x, x_{\mu})$ and is decreasing in $(x_\mu, \bar x+1)$. Hence $\tilde w$ must be uniquely given by the formula (\ref{physical-sol}). \qed

\bibliographystyle{plain}

\end{document}